\documentclass[]{elsarticle}

\journal{Journal of Pure and Applied Algebra}
\usepackage{color,amssymb,amsmath,amsthm,tikz,pb-diagram,mathtools,natbib}
\usepackage{mathrsfs}

\usepackage{enumerate}
\usepackage[shortlabels]{enumitem}

\usepackage{hyperref}
\hypersetup{
 unicode=false, 
 pdftoolbar=true, 
 pdfmenubar=true, 
 pdffitwindow=false, 
 pdfstartview={FitH}, 
 pdftitle={comments32}, 
 pdfauthor={JGG,UH,TK}, 
 pdfsubject={points}, 
 pdfcreator={TeXShop}, 
 pdfproducer={Producer}, 
 pdfkeywords={keywords}, 
 pdfnewwindow=true, 
 colorlinks=true, 
 linkcolor=altncolor, 
 citecolor=green, 
 filecolor=blue, 
 urlcolor=cyan 
}

\journal{}

 \newcommand{\twoheaddownarrow}{{\rlap{\rlap{$\ $}\raise .25ex\hbox{$\downarrow$}}\raise-.25ex\hbox{$\downarrow$}}}
 \newcommand{\twoheaduparrow}{{\rlap{\rlap{$\ $}\raise .25ex\hbox{$\uparrow$}}\raise-.25ex\hbox{$\uparrow$}}}

\newcommand{\tbigwedge}{\mathop{\textstyle \bigwedge }}
\newcommand{\tbigcup}{\mathop{\textstyle \bigcup }}
\newcommand{\tbigvee}{\mathop{\textstyle \bigvee }}

\def\LR{{\mathfrak{L}(\R)}}

\def\N{{\mathbb{N}}}

\def\R{{\mathbb{R}}}
\def\Q{{\mathbb{Q}}}

\def\Frm{{\mathsf{Frm}}}
\def\Top{{\mathsf{Top}}}
\def\Pt{{\Sigma}}
\def\O{{\mathcal O}}

 \newcommand{\gidoia}{\textsf{\scriptsize ---}}
 \newcommand{\cbelow}{{\prec\hskip-5pt \prec\hskip2pt}}

\newtheorem{theorem}{Theorem}[section]
\newtheorem{proposition}[theorem]{Proposition}
\newtheorem{lemma}[theorem]{Lemma}
\newtheorem{corollary}[theorem]{Corollary}
\newtheorem{fact}[theorem]{Fact}
\theoremstyle{definition}
\newtheorem{definition}[theorem]{Definition}
\newtheorem{example}[theorem]{Example}
\newtheorem{examples}[theorem]{Examples}

\theoremstyle{remark}
\newtheorem{remark}[theorem]{Remark}
\newtheorem{remarks}[theorem]{Remarks}

\makeatletter
\newcommand*{\@old@slash}{}\let\@old@slash\slash
\def\slash{\relax\ifmmode\delimiter"502F30E\mathopen{}\else\@old@slash\fi}
\makeatother

\DeclareMathOperator{\In}{I}
\DeclareMathOperator{\Su}{S}

\begin{document}

\begin{frontmatter}

\title{On the universal completion\\ of pointfree function spaces
}

\author{Imanol~Mozo~Carollo}
 \ead{imanol.mozo@ehu.eus}


\address{
Department of Applied Economics I, Faculty of Economics and Bussiness, \\
University of the Basque Country UPV-EHU, \\Plaza O\~nati 1, 20018 Donostia-San Sebasti\'an, Spain}


 \begin{abstract}
This paper approaches the construction of the universal completion of the Riesz space $\mathrm{C}(L)$ of continuous real functions on a completely regular frame $L$ in two different ways. Firstly as the space of continuous real functions on the Booleanization of $L$. Secondly as the space of nearly finite Hausdorff continuous functions on $L$. The former has no counterpart in the classical theory, as the Booleanization  of a spatial frame is not spatial in general, and it offers a lucid way of representing the universal completion as a space of continuous real functions. As a corollary we obtain that $\mathrm{C}(L)$ and $\mathrm{C}(M)$ have isomorphic universal completions if and only if the Booleanization of $L$ and $M$ are isomorphic and we characterize frames $L$ such that $\mathrm{C}(L)$ is universally complete as almost Boolean frames. The application of this last result to the classical case $\mathrm{C}(X)$ of the space of continuous real functions on a topological space $X$ characterizes those spaces $X$ for which $\mathrm{C}(X)$ is universally complete. Finally, we present a pointfree version of the Maeda-Ogasawara-Vulikh representation theorem and use it to represent the universal completion of an Archimedean Riesz space with weak unit as a space of continuous real functions on a Boolean frame.
 \end{abstract}


\end{frontmatter}

 \section*{Introduction} Pointfree topology replaces classical spaces by an abstraction of their lattices of open subsets, namely, \emph{frames} (or \emph{locales}). This approach to topology is not merely a restatement of the classical theory,  but an actual generalization in which spaces ``without enough points'' are also allowed. Our main goal with this paper is to take advantage of this feature in order to present a construction of the universal completion of the Riesz space $\mathrm{C}(X)$ of continuous real functions on a space $X$ in the most lucid and simple way.

We begin by recalling that a Riesz space $R$  is called: \emph{Dedekind complete} if each non-void subset $A\subseteq R$ that is bounded from above has a supremum in $R$; \emph{laterally complete} if each non-void set $A\subseteq R^+=\{f\in R\mid 0\leq f\}$ of disjoint elements, that is, such that $g\wedge h=0$ for each $g,h\in A$, has a supremum in $R$; and \emph{universally complete} if it is both Dedekind complete and laterally complete.
The \emph{universal completion}  of a Riesz space $R$ is  a pair $(\mathcal{U}(R),\mu)$ where $\mathcal{U}(R)$ is a universally complete Riesz space and  $\mu\colon R\to\mathcal{U}(R)$ is a Riesz space embedding such that for every $f\in \mathcal{U}(R)^+$,
\[
f=\tbigvee\{\mu(g)\colon g\in R,\, 0\leq \mu(g)\leq f\}.
\]
The existence of the universal completion of an Archimedean Riesz space traces back to \cite{MO42} where  Maeda and Ogasawara presented a description of the universal completion based on a representation of Riesz spaces as spaces of continuous extended real functions (see \cite{LZ71}). Uniqueness is due to Vulikh  \cite{KVP50}. For the particular case of the space $\mathrm{C}(X)$ of continuous real functions on a completely regular space $X$, van der Walt presented in \cite{W18} a direct construction of its universal completion in terms of nearly finite lower semicontinuous functions on the same space. 
 
In this paper, we address the universal completion of the space of continuous real functions from the point of view of pointfree topology. In this setting, frame homomorphisms $f\colon L\to M$ represent continuous functions between the frames  $L$ and $M$, thought of as generalized spaces. For a frame $L$, the \emph{continuous real functions} on $L$ are all the frame homomorphisms from the \emph{frame of reals} $\LR$ into $L$. They naturally form a Riesz space $\mathrm{C}(L)$ \cite{BB97}. The correspondence $L\mapsto \mathrm{C}(L)$ extends that for spaces: if $L=\mathcal{O}X$ (the frame of open sets of a space $X$) then the classical function space $\mathrm{C}(X)$ is naturally isomorphic to $\mathrm{C}(L)$ (see \cite{BB97}).

After studying in Section 2 the universal completeness of $\mathrm{C}(L)$ when $L$ is Boolean, we prove that the natural embedding of $\mathrm{C}(L)$ into $\mathrm{C}(\mathfrak{B}(L))$ constitutes its universal completion (Section 3), where $\mathfrak{B}(L)$ denotes the Boolenization of $L$. The Booleanization of a frame is the natural generalization of the lattice of regular open subsets of a classical space. While $\mathfrak{B}(L)$ is always a frame and, consequently, a perfectly admissible ``space'' in the sense of pointfree topology, this constructions has no analogue in the classical theory, as the Booleanization of a frame $L$ may not be spatial even in the case of a spatial frame $L$. The broader class of spaces that pointfree topology offers is precisely the reason why our approach provides a more direct construction and sheds new light on the problem.

In Section 4 we present an alternative description of the completion in terms of a certain kind of interval valued functions, namely, \emph{nearly finite Hausdorff continuous functions}.  Then we make use of our representations in order to characterize those completely regular frames $L$ for which $\mathrm{C}(L)$  is universally complete, namely, almost Boolean frames, and to show that $\mathrm{C}(L)$ and $\mathrm{C}(M)$ have isomorphic universal completions iff $\mathfrak{B}(L)$ and $\mathfrak{B}(M)$ are isomorphic (Section 5). We use this construction later in Section 6 in the application of these ideas to the classical case $\mathrm{C}(X)$, providing yet another representation in the spatial setting. The construction of van der Walt follows easily from our approach. Furthermore, we obtain a direct proof that avoids representation theory for the following result of van der Walt \cite[Corollary 39]{W18}: $\mathrm{C}(X)$ and $\mathrm{C}(Y)$ have Riesz isomorphic universal completions if and only if $\mathfrak{B}(\O X)\simeq\mathfrak{B}(\O Y)$.

Finally, in Section 7 we give account of the classical Maeda-Ogasawara-Vulikh representation theorem for Archimedean Riesz spaces with weak unit from a pointfree point of view. Specifically, we show that for an Archimedean Riesz space $R$ with weak unit, $R$ embeds into the  space $\mathrm{C}(\mathcal{B}(R))$ of continuous real functions on the frame $\mathcal{B}(R)$ of bands of $R$. In addition, this embedding is precisely the universal completion of $R$. It is worth pointing out that this representation not only does not need extended real valued functions, but also avoids the construction of the Stone space of $\mathcal{B}(R)$, leading to a simpler construction.
 
 \section{Background}
 
\subsection{Notation}
For any subset $A$ of a partially ordered set $(P,\le)$, we will denote by $
 {\tbigvee^P} A$ (resp. $
 {\tbigwedge^P} A$) the supremum (resp. infimum) of $A$ in $P$ in case it exists (we shall omit the superscript if it is clear from the context).

\subsection{Frames}

A \emph{frame}\index{frame} (or \emph{locale}\index{locale}) $L$ is a complete lattice such that
 \begin{equation}
 a\wedge \tbigvee B=\tbigvee \{a\wedge b\mid b\in B\}
 \end{equation}
 for all $a\in L$ and $B\subseteq L$; equivalently, it is a complete Heyting algebra with Heyting operation $\rightarrow$ satisfying the standard equivalence $a\wedge b\le c$ if and only if $a\le b\rightarrow c$. The \emph{pseudocomplement} of an $a\in L$ is the element
 \[
 a^{\ast}=a\rightarrow 0=\tbigvee \{b\in L\mid a\wedge b=0\}.
 \]
An element $a\in L$ is \emph{complemented} if $a\vee a^\ast=1$ and \emph{dense} if $a^\ast=0$ (equivalently, if $a^{\ast\ast}=1$). A \emph{frame homomorphism}\index{frame homomorphism} is a map $h\colon L\to M$ between frames which preserves finitary meets (including the top element 1) and arbitrary joins (including the bottom element 0). Then $\Frm$ is the corresponding category of frames and their homomorphisms.

 The most typical example of a frame is the lattice $\O X$ of open subsets of a topological space $X$. The correspondence $X\mapsto \O X$ is clearly functorial, and consequently we have a contravariant functor $\O \colon \Top\to\Frm$ where $\Top$ denotes the category of topological spaces and continuous maps. There is also a functor in the opposite direction, the \emph{spectrum functor}\index{spectrum functor}\index{functor!spectrum $\sim$}\index{functor!$\Pt$} $\Sigma\colon\Frm\to\Top$ which assigns to each frame $L$ its spectrum $\Sigma L$, the space of all homomorphisms $\xi\colon L\to \{0,1\}$ with open sets $\Sigma_a=\{\xi\in\Sigma L\mid \xi(a)=1\}$ for any $a\in L$, and to each frame homomorphism $h\colon L\to M$ the continuous map $\Sigma h\colon\Sigma M\to\Sigma L$ such that $\Sigma h(\xi)=\xi h$. The spectrum functor is right adjoint to $\mathcal O$, with adjunction maps $\eta_L\colon L\to \O \Sigma L, \;\eta_L(a)=\Sigma_a$ and $\epsilon_X\colon X\to \Sigma
 \O X, \; \epsilon_X(x)=\hat{x}, \; \hat{x}(U)=1\mbox{ if and only if }x\in U$ (the former is the \emph{spatial reflection} of the frame $L$). A frame is said to be \emph{spatial}\index{spatial frame}\index{frame!spatial $\sim$} if it is isomorphic to the frame of open sets of a space.

For general notions and results concerning frames we refer to Johnstone \cite{PJ82} or the recent
Picado-Pultr \cite{PP12}. The particular notions we will need are the following:
a frame $L$ is:

 \begin{enumerate}[-]
%


 \item \emph{regular} if $a = \tbigvee\{b\in L\mid b \prec a\}$ for every $a\in L$, where $b\prec a$ ($b$ is \emph{rather below}\index{rather below relation}\index{relation!rather below $\sim$} $a$) means that $b^\ast\vee a=1$;

 \item \emph{completely regular} if $a = \tbigvee\{b\in L\mid b{\prec\hskip-2pt \prec\hskip2pt}a\}$ for each $a\in L$, where $b{\prec\hskip-2pt \prec\hskip2pt}a$ ($b$ is \emph{completely below}\index{completely below relation}\index{relation!completely below $\sim$} $a$) means that there is $\{c_{r}\mid r\in \mathbb{Q}\cap[0,1]\}\subseteq L$ such that $b\le c_{0}$, $c_{1}\le a$ and $c_r\prec c_s$ whenever $r<s$;

\item \emph{extremally disconnected} if $a^\ast\vee a^{\ast\ast}=1$ for every $a\in L$.
 \end{enumerate}

 \subsection{Continuous real functions}
 
 One of the main differences between pointfree topology and classical topology is that the category of
frames is algebraic. In consequence, we can present frames  by specifying generators and relations over those generators in terms of finite meets and arbitrary joins. This fact was used by Joyal in order to introduce the \emph{frame $\mathfrak{L}(\R)$ of reals}, the pointfree counterpart of the real line \cite{J73}. This frame was further studied by Banaschewski in \cite{BB97}, with a special emphasis on the rings of continuous real functions.

Here it will be useful to adopt the equivalent description of $\mathfrak{L}(\R)$ used in \cite{GKP09} given by generators $(p,{\gidoia})$ and $({\gidoia},p)$, $p\in\Q$ subject to the defining relations
\begin{enumerate}[xxxxx] 
\item[{\rm(r1)}] $(p,{\gidoia})\wedge ({\gidoia},q)=0$ whenever $p\ge q$,
\item[{\rm(r2)}] $(p,{\gidoia})\vee ({\gidoia},q)=1$ whenever $p<q$,
\item[{\rm(r3)}] $(p,{\gidoia})= \bigvee_{q>p}(q,{\gidoia})$, for every $p\in\Q$,
\item[{\rm(r4)}] $({\gidoia},p)= \bigvee_{q<p}({\gidoia},q)$, for every $p\in \Q$,
\item[{\rm(r5)}] $ \bigvee_{p\in \mathbb{Q}}(p,{\gidoia})=1$,
\item[{\rm(r6)}] $\bigvee_{p\in \mathbb{Q}}({\gidoia},p)=1$.
\end{enumerate}
The meet $(p,\gidoia)\wedge (\gidoia,q)$ is simply denoted by $(p,q)$. By dropping relations (r5) and (r6) in the description of $\mathfrak{L}(\mathbb{R})$ above, we have the corresponding \emph{frame of extended reals} $\mathfrak{L}\big(\overline{\mathbb{R}}\big)$ \cite{BGP14}. 

For any frame $L$, a \emph{continuous real function} \cite{BB97} (resp. \emph{extended continuous real function}\index{continuous!extended real function}\cite{BGP14}) on a frame $L$ is a frame homomorphism $f\colon {\mathfrak{L}(\mathbb{R})}\rightarrow L$ (resp. $f\colon {\mathfrak{L}\big(\overline{\mathbb{R}}\big)}\rightarrow L$). We denote by $\mathrm{C}(L)$ (resp. $\overline{\mathrm{C}}(L)$) the collection of all (resp. extended) continuous real functions on $L$. Considering frame homomorphisms from $\mathfrak{L}(\R)$ (rexp. $\mathfrak{L}\big(\overline{\R}\big)$) into a general frame $L$ as continuous (resp. extended) real functions of $L$ provides a natural extension of the classical notion as contiuous (resp. extended) real functions on a space $X$ can be represented by homomorphisms $\mathfrak{L}(\R)\to \O X$ (resp. $\mathfrak{L}\big(\overline{\R}\big)\to\O X$) (see \cite{BB97} and \cite{BGP14} for details). There is a basic homomorphism ${\varrho}\colon {\mathfrak{L}\big(\overline{{{\mathbb{R}}}}\big)}\to {{\mathfrak{L}({{\mathbb{R}}})}}$, determined on generators by
 \[
 \varrho(p,\gidoia)=(p,\gidoia)\quad\text{and}\quad\varrho(\gidoia,q)=(\gidoia,q)
 \]
 for each $p,q\in\Q$. The functions $f\in \overline{\mathrm{C}}(L)$ that factor through the canonical homomorphism $\varrho$ are just those that turn the defining relations  (r5) and (r6) into identities in $L$. Accordingly, we will still denote by $\mathrm{C}(L)$  the class inside $\overline{\mathrm{C}}(L)$ of functions that turn (r5) and (r6) into identities in $L$.

$\mathrm{C}(L)$ and $\overline{\mathrm{C}}(L)$ are partially ordered by
\begin{equation}\label{fleg}
\begin{aligned}
f\le g& \iff f(p,\gidoia)\le g(p,{\gidoia})\quad\mbox{for all }p\in{{\mathbb{Q}}}\\
&\iff g({\gidoia},q)\le f({\gidoia},q)\quad\mbox{for all }q\in{{\mathbb{Q}}}.
\end{aligned}
\end{equation}
It is well known that $\mathrm{C}(L)$ is not Dedekind complete in general. In fact $\mathrm{C}(L)$ is Dedekind complete if and only if $L$ is extremally disconnected \cite{BH03}.

\begin{remark}\label{funtziobikoitza}
For any $f\in\mathrm{C}(L)$ and $p\in\Q$,
\[
f(p,\gidoia)=\tbigvee_{r>p}f(r,\gidoia)^{\ast\ast}.
\]
Obviously $f(p,\gidoia)\leq\tbigvee_{r>p}f(r,\gidoia)^{\ast\ast}$ as $f(r,\gidoia)^{\ast\ast}\geq f(r,\gidoia)$ for any $r\in\Q$. In order to check that $f(p,\gidoia)\geq\tbigvee_{r>p}f(r,\gidoia)^{\ast\ast}$, note that $f(r,\gidoia)^{\ast\ast}\leq f(p,\gidoia)$ for any $r>p$. Indeed
\[
\begin{aligned}
f(r,\gidoia)^{\ast\ast}&=(f(p,\gidoia)\vee f(\gidoia,r))\wedge f(r,\gidoia)^{\ast\ast}\\
&=(f(p,\gidoia)\wedge f(r,\gidoia)^{\ast\ast})\vee(f(\gidoia,r)\wedge f(r,\gidoia)^{\ast\ast})\\
&=f(p,\gidoia)\wedge f(r,\gidoia)^{\ast\ast}
\end{aligned}
\]
since $f(\gidoia,r)\wedge f(r,\gidoia)^{\ast\ast}=0$ iff $f(\gidoia, r)\wedge f(r,\gidoia)=0$.
\end{remark}

\begin{example}\label{constantfunction}
For each $r\in \mathbb{Q}$, the \emph{constant function} $\boldsymbol{r}$ determined by $r$ is defined by
\[\boldsymbol{r}(p,\gidoia)= \begin{cases}0&\text{if }p\ge r,\\ 1&\text{if }p< r, \end{cases}\quad\text{ and }\quad\boldsymbol{r}(\gidoia,q)= \begin{cases}1&\text{if }q>r,\\ 0&\text{if }q\le r, \end{cases}\]
for every $p,q\in{{\mathbb{Q}}}$.
\end{example}
 
  \subsection{Scales.}
There is a useful way of specifying (extended) continuous real functions on a frame $L$ with the help of the so called (extended) scales (\cite[Section~4]{GKP09}).
An \emph{extended scale} in $L$ is a map $\sigma\colon {{\mathbb{Q}}}\rightarrow L$ such that $\sigma(q)\prec\sigma(p)$ whenever $p<q$. An extended scale is a \emph{scale} if
\[
\tbigvee_{p\in \mathbb{Q}}\sigma(p)=1=\tbigvee_{p\in \mathbb{Q}}\sigma(p)^\ast.
\]
For each extended scale $\sigma$ in $L$, the formulas
\begin{equation}\label{fscales}
f(p,{\gidoia})=\tbigvee_{r>p}\sigma(r)\quad \text{and}\quad f({\gidoia},q)= \tbigvee_{r<q}{\sigma(r)}^\ast,
\quad p,q\in\Q,
\end{equation}
determine an $f\in \overline{\mathrm{C}}(L)$; then,
$f\in\mathrm{C}(L)$ if and only if $\sigma$ is a scale.

An extended scale $\sigma$ is necessarily antitone. However, a map $\sigma\colon\Q\to L$ such that $\sigma(p)$ is  complemented for every $p\in\Q$ is an
extended scale if and only if it is antitone.

 \subsection{Algebraic operations on $\mathrm{C}(L)$}
 
 The operations on the algebra $\mathrm{C}(L)$ are determined by the operations of $\Q$ as an ordered vector lattice as follows (see \cite{BB97} for more details):
 \begin{enumerate}[\rm(1)]
 \item[\rm(1)] For $\diamond=+,\cdot,\wedge,\vee$:
 $$(f\diamond g)(p,q)=\tbigvee\{f(r,s)\wedge g(t,u)\mid \langle r,s\rangle\diamond\langle t,u\rangle\subseteq \langle p,q\rangle\}$$
where $\langle\cdot,\cdot\rangle$ stands for open interval in $\mathbb{Q}$ and the inclusion on the right means
that $x\diamond y\in\langle p,q\rangle$ whenever
$x\in \langle r,s\rangle$ and $y\in\langle t,u\rangle$.
 \item[\rm(2)] $({\scriptstyle-}f)(p,q)=f(-q,-p)$.
 \item[\rm(3)] For each $r\in \Q$, the nullary operation $\boldsymbol{r}$ is defined as in Example \ref{constantfunction} above.
 \item[\rm(4)] For each $0<\lambda\in \mathbb{Q}$, $(\lambda\cdot f)(p,q)=f\left(\tfrac{p}{\lambda},\tfrac{q}{\lambda}\right)$.
 \end{enumerate}
 These operations satisfy all the identities which hold for their counterparts in $\mathbb{Q}$ and hence they determine an
 Archimedean Riesz space structure with weak unit $\boldsymbol{1}$ in $\mathrm{C}(L)$ (see \cite{M90}). 
 The following formulas will be useful \cite{GP11}:
  \begin{enumerate}[\rm(1)]
  \item $(f\wedge g)(p,\gidoia)=f(p,\gidoia)\wedge g(p,\gidoia)$ and $(f\wedge g)(\gidoia,q)=f(\gidoia,q)\vee g(\gidoia,q)$;
  \item $(f\vee g)(p,\gidoia)=f(p,\gidoia)\vee g(p,\gidoia)$ and $(f\vee g)(\gidoia,q)=f(\gidoia,q)\wedge g(\gidoia,q)$;
  \item $(f+g)(p,\gidoia)=\tbigvee_{r\in\Q}f(r,\gidoia)\wedge g(p-r,\gidoia)$ and
  \newline $(f+g)(\gidoia,q)=\tbigvee_{s\in\Q}f(\gidoia, s)\wedge g(\gidoia, t-s)$;
 \item For each $0<\lambda\in\Q$, $(\lambda\cdot f)(p,\gidoia)=f\left(\tfrac{p}{\lambda},\gidoia\right)$ and $f\left(\gidoia,\tfrac{q}{\lambda}\right)$.
  \end{enumerate}

 \section{Continuous real functions on Boolean frames}
 
 Boolean frames, that is, frames in which every element is complemented, are the natural counterpart of discrete spaces in pointfree topology \cite{PP17}. This is a proper extension of the classical notion, as not all Boolean frames are spatial. In fact, a complete Boolean algebra is a spatial frame if and only if it is atomic \cite{PP12}. In this section we establish the main results that we will need later about spaces of continuous real functions on Boolean frames. 
 
Recall that by a \emph{discrete} $\{y_i\}_{i\in I}$ in $L$ it is meant a collection for which there is a cover $C$ of $L$, that is, $\tbigvee C=1$, such that for each $c\in C$, $c\wedge y_i=0$ for all $i$ with possibly one exception. By a \emph{co-discrete} $\{{x_i}\}_{i\in I}$ it is meant a collection for which there is a cover $C$ such that for each $c\in C$, $c\le x_i$ for all $i$ with possibly one exception. A collection $\{y_i\}_{i\in I}$ is discrete if and only if $\{y_i^*\}_{i\in I}$ is co-discrete and for any co-discrete system $\{x_i\}_{i\in I}$ and any $y\in L$, $y\vee \tbigwedge_{i\in I} x_i=\tbigwedge_{i\in I}(y\vee x_i)$ \cite{AP84}. Recall also that $\mathrm{C}(L)^+=\{f\in\mathrm{C}(L)\mid \boldsymbol{0}\leq f\}$.
 
 \begin{proposition}\label{discretejoin}
 Let $L$ be frame and $\{f_i\}_{i\in I}\subseteq \mathrm{C}(L)^+$ be such that $\{f_i(0,\gidoia)\}_{i\in I}$ is discrete. Then $\tbigvee f_i$ exists in $\mathrm{C}(L)$.
 \end{proposition}
 
\begin{proof}
 Let $\sigma\colon\Q\to L$ be the map given by $\sigma(p)=\tbigvee_{i\in I}f_i(p,\gidoia)$ for each $p$ in $\Q$. If $0\leq p<q$ in $\Q$, one has
 \[
 \sigma(p)\vee\sigma(q)^\ast=\tbigvee_{i\in\Q}f_i(p,\gidoia)\vee\tbigwedge_{i\in I}f(q,\gidoia)^\ast=\tbigwedge_{i\in I}\left(\tbigvee_{i\in\Q}f_i(p,\gidoia)\vee f(q,\gidoia)^\ast\right)=1,
 \]
 since $\{f_i(q,\gidoia)^\ast\}_{i\in I}$ is co-discrete as $\{f_i(0,\gidoia)^\ast\}_{i\in I}$ is co-discrete and $f_i(0,\gidoia)^\ast\leq f_i(q,\gidoia)^\ast$ for each $i\in I$. If $p<0$ then $\sigma(p)=1$. Hence $\sigma(q)\prec\sigma(p)$ whenever $p<q$. Now, let $C$ be the cover that witnesses  the co-discreteness of $\{f_i(0,\gidoia)^\ast\}_{i\in I}$ and $c\in C$. We have
 \begin{enumerate}[(a)]
 \item If $c\leq f_i(0,\gidoia)^\ast$ for all $i\in I$, then
 \[
 \tbigvee_{p\in\Q}\sigma(p)^\ast=\tbigvee_{p\in\Q}\tbigwedge_{i\in I}f_i(p,\gidoia)^\ast \geq\tbigvee_{p\geq 0}\tbigwedge_{i\in I}f_i(0,\gidoia)^\ast\geq c.
 \]
 \item If $c\not\leq f_{i_0}(0,\gidoia)^\ast$ for some $i_0\in I$, note that $c\wedge f_i(p,\gidoia)^\ast=c$ for every $i\neq i_0$ as $f_i(0,\gidoia)\leq f_i(p,\gidoia)^\ast$. Accordingly,
 \[
 \begin{aligned}
 c\wedge\tbigvee_{p\in\Q}\sigma(p)^\ast&=c\wedge\tbigvee_{p\in\Q}\tbigwedge_{i\in I}f_i(p,\gidoia)^\ast =\tbigvee_{p\in\Q}\tbigwedge_{i\in I}(c\wedge f_i(p,\gidoia)^\ast)\\
 &\geq\tbigvee_{p\geq 0}\tbigwedge_{i\in I}(c\wedge f_i(p,\gidoia)^\ast)=\tbigvee_{p\geq 0}(c\wedge f_{i_0}(p,\gidoia)^\ast)\\
 &=c\wedge\tbigvee_{p\geq 0}f_{i_0}(p,\gidoia)^\ast=c,
 \end{aligned}
 \]
 Therefore $c\leq\tbigvee_{p\in\Q}\sigma(p)^\ast$.
 \end{enumerate}
 Hence $\tbigvee_{p\in\Q}\sigma(p)^\ast\geq\tbigvee_{c\in C}c=1$. It is obvious that $\tbigvee_{p\in\Q}\sigma(p)=1$. Consequently $\sigma$ is a scale on $L$ that obviously determines the supremum of $\{f_i\}_{i\in I}$.
\end{proof}

 \begin{corollary}\label{Booleanunivcomp}
 For any Boolean frame $L$, $\mathrm{C}(L)$ is universally complete.
 \end{corollary}
 
 \begin{proof}
 First note that $\mathrm{C}(L)$  is Dedekind complete since $L$ is extremally disconnected. In order to check that $\mathrm{C}(L)$ is laterally complete, simply note $\{f_i\}_{i\in I}$ is pairwise disjoint iff $\{f_i(0,\gidoia)\}_{i\in I}$ is pairwise disjoint and pairwise disjointness is equivalent to discreteness in Boolean frames. Thus  $\tbigvee_{i\in I}f_i$ exists in $\mathrm{C}(L)$ by Proposition \ref{discretejoin}
 \end{proof}

The following is certainly folklore, but as far as the author knows it has not been published. For the sake of completeness, we include the proof here, as this result will be useful later in order to study when the universal completions of $\mathrm{C}(L)$ and $\mathrm{C}(M)$ are isomorphic.
 
 \begin{proposition}\label{isoBooleanpart}
 Let $L$ and $M$ be Boolean frames. If $\mathrm{C}(L)$ and $\mathrm{C}(M)$ are isomorphic, so are $L$ and $M$.
 \end{proposition}
 
 \begin{proof}
  First, if $\Phi\colon \mathrm{C}(L)\to \mathrm{C}(M)$ is a lattice isomorphism such that $\Phi(\boldsymbol{0})=\boldsymbol{0}$, we shall check that $\varphi\colon L\to M$, where $\varphi(a)=\Phi(\chi_a)(0,\gidoia)$, is also an isomorphism. For any $a,b\in L$ one has
 \[
 \begin{aligned}
 \varphi(a)\wedge\varphi(b)&=\Phi(\chi_a)(0,\gidoia)\wedge\Phi(\chi_b)(0,\gidoia)\\
 &=(\Phi(\chi_a)\wedge\Phi(\chi_b))(0,\gidoia)\\
 &=\Phi(\chi_a\wedge\chi_b)(0,\gidoia)\\
 &=\Phi(\chi_{a\wedge b})(0,\gidoia)=\varphi(a\wedge b).
 \end{aligned}
 \]
 One can check $\varphi(a)\vee\varphi(b)=\varphi(a\vee b)$ dually.
 
We shall check that, for any $f\in\mathrm{C}(L)^+$, if $f(0,\gidoia)=1$ then $\Phi(f)(0,\gidoia)=1$. Let $\Phi(f)(0,\gidoia)=a$. Then we have $\Phi^{-1}(\chi_{a^c})\wedge f=\boldsymbol{0}$, as $\chi_{a^c}\wedge\Phi(f)=\boldsymbol{0}$. Therefore $0=(\Phi^{-1}(\chi_{a^c})\wedge f)(0,\gidoia)=\Phi^{-1}(\chi_{a^c})(0,\gidoia)\wedge f(0,\gidoia)=a^c\wedge 1$. Hence $a=1$. In particular one has $\varphi(1)=\Phi(\chi_1)=1$.  As we obviously have $\varphi(0)=0$, $\varphi$ is a Boolean algebra homomorphism.
 
 In order to check that $\varphi$ is onto, let $a\in M$ and $b=\Phi^{-1}(\chi_a)(0,\gidoia)\in L$. Then we have $\chi_{b^c}\wedge\Phi^{-1}(\chi_a)=\boldsymbol{0}$ and $(\chi_{b^c}\vee\Phi^{-1}(\chi_a))(0,\gidoia)=1$. Consequently $\varphi(b^c)\wedge a=0$ and $\varphi(b^c)\vee a=1$. Hence, $\varphi(b)=a$. Thus $\varphi$ is onto. Now let $a\neq b\in L$. Then one has that, say, $a^c\wedge b\neq 0$. Hence, $\varphi(a^c)\wedge\varphi(b)\neq 0$ as $\chi_{a^c}\wedge\chi_b\neq\boldsymbol{0}$. Therefore $\varphi(a)\neq\varphi(b)$. Consequently, $\varphi$ is injective.

 Finally, for the general case where the isomorphism $\Phi\colon\mathrm{C}(L)\to\mathrm{C}(M)$ not necessarily maps $\boldsymbol{0}$ into $\boldsymbol{0}$, simply note that $\overline{\Phi}\colon\mathrm{C}(L)\to\mathrm{C}(M)$ given by $\overline{\Phi}(f)=\Phi(f)-\Phi(\boldsymbol{0})$ is also an isomorphism. 
 \end{proof}
 
 For the question whether the last proposition holds for a broader class of frames, see Remarks \ref{remarkAlmostBoolean} (1) below.
 
 \section{The universal completion of $\mathrm{C}(L)$}
 
 In this section we establish the main result of the paper. We show that the universal completion of the space of continuous real functions on any completely regular frame $L$ is isomorphic to the lattice of all continuous real functions on the \emph{Booleanization} of $L$:
\[
\mathfrak{B}(L)=\{a\in L\mid a=a^{\ast\ast}\}.
\]
This is
the complete Boolean algebra of all regular elements $a=a^{**}$ of $L$. Meets in $\mathfrak{B}(L)$ coincide with meets in $L$ while joins are given by
\[
\tbigvee^{\mathfrak{B}(L)}A=\left(\tbigvee^L A\right)^{\ast\ast}
\]
for each $A\subseteq L$. There is an associated frame homomorphism $\beta\colon L\to \mathfrak{B}(L)$ that maps each element $a$ to its double pseudocomplement $a^{\ast\ast}$ \cite{BP96}.

 \begin{proposition}
 The map $\Upsilon\colon\mathrm{C}(L)\to\mathrm{C}(\mathfrak{B}(L))$ given by $f\mapsto \beta\cdot f$ is a Reisz embedding.
 \end{proposition}
 
 \begin{proof}
First we check that $\Upsilon$ is a Riesz homomorphism. Obviously, $\Upsilon(\boldsymbol{0})=\boldsymbol{0}$. For any $0<\lambda\in\Q$ and $f\in\mathrm{C}(L)$, we have \[
 \begin{aligned}
 (\lambda\cdot\Upsilon(f))(p,\gidoia)&=\Upsilon(f)\left(\tfrac{p}{\lambda},\gidoia\right)=f\left(\tfrac{p}{\lambda},\gidoia\right)^{\ast\ast}=\Upsilon(\lambda\cdot f)(p,\gidoia)
 \end{aligned}
 \]
 for all $p\in\Q$. Consequently, $\lambda\cdot \Upsilon(f)=\Upsilon(\lambda\cdot f)$. Further, for any $f,g\in\mathrm{C}(L)$, we have
 \[
 \begin{aligned}
 (\Upsilon(f)+\Upsilon(g))(p,\gidoia)&=\tbigvee^{\mathfrak{B}(L)}_{r\in\Q}(f(r,\gidoia)^{\ast\ast}\wedge g(p-r,\gidoia)^{\ast\ast})\\
 &=\left(\tbigvee^L_{r\in\Q}(f(r,\gidoia)^{\ast\ast}\wedge g(p-r,\gidoia)^{\ast\ast})\right)^{\ast\ast}\\
 &=\left(\tbigvee^L_{r\in\Q}(f(r,\gidoia)\wedge g(p-r,\gidoia))^{\ast\ast}\right)^{\ast\ast}\\
 &=\left(\tbigwedge^L_{r\in\Q}(f(r,\gidoia)\wedge g(p-r,\gidoia))^{\ast\ast\ast}\right)^{\ast}\\
 &=\left(\tbigvee^L_{r\in\Q}(f(r,\gidoia)\wedge g(p-r,\gidoia))\right)^{\ast\ast}\\
 &=(f+g)(p,\gidoia)^{\ast\ast}=\Upsilon(f+g)(p,\gidoia)
 \end{aligned}
 \]
 for all $p\in\Q$. Consequently, $\Upsilon(f)+\Upsilon(g)=\Upsilon(f+g)$. For suprema, let $f,g\in\mathrm{C}(L)$ and we have
 \[
 \begin{aligned}
 (\Upsilon(f)\vee\Upsilon(g))(p,\gidoia)&=\Upsilon(f)(p,\gidoia)\vee^{\mathfrak{B}(L)} \Upsilon(g)(p,\gidoia)\\
 &=(f(p,\gidoia)^{\ast\ast}\vee g(p,\gidoia)^{\ast\ast})^{\ast\ast}\\
 &=(f(p,\gidoia)\vee g(p,\gidoia))^{\ast\ast}=\Upsilon(f\vee g)(p,\gidoia)
 \end{aligned}
 \]
 for all $p\in\Q$. Hence $\Upsilon(f)\vee\Upsilon(g)=\Upsilon(f\vee g)$. One can check dually that $\Upsilon(f)\wedge\Upsilon(g)=\Upsilon(f\wedge g)$. Finally, injectivity follows easily from Remark \ref{funtziobikoitza}.

 \end{proof}
 
 \begin{theorem}\label{firstcompletion}
For a completely regular frame $L$, $\Upsilon\colon\mathrm{C}(L)\to\mathrm{C}(\mathfrak{B}(L))$ is the universal completion of $\mathrm{C}(L)$.
 \end{theorem}
 
 \begin{proof}
 Let $\boldsymbol{0}\leq h\in\mathrm{C}(\mathfrak{B}(L))$ and $\mathcal{F}=\{f\in\Upsilon(\mathrm{C}(L))\mid \boldsymbol{0}\leq f\leq h\}$. 
 Since $\mathrm{C}(\mathfrak{B}(L))$ is Dedekind complete, the supremum $f_\vee=\tbigvee^{\mathrm{C}(\mathfrak{B}(L))}\mathcal{F}$ exists. We shall prove that $f_\vee=h$.
 
 For this purpose, fix $0\leq q\in\Q$ and an $a\in L$ such that $a\cbelow h(q,\gidoia)$.  Then there exists a family $\{c_r\mid r\in\Q\cap[0,1]\}$ such that $a\leq c_0, c_1\leq h(q,\gidoia)$ and $c_r\prec c_s$ whenever $r<s$. Let $\psi_q\colon \Q\cap[0, q]\to[0,1]$ be a dual-order isomorphism and define the mapping $\sigma_{q,a}\colon\Q\to L$ by
 \[
 \sigma_{q,a}(p)=\begin{cases} 
 1&\text{if }p<0\\
 c_{\psi(p)}&\text{if }0\leq p< q \\
 0&\text{if } q\leq p.
 \end{cases}
 \]
 It is straighforward to check thath $\sigma_{q,a}$ is a scale on $L$ and consequently it determines a continuous real function $f_{q,a}\in\mathrm{C}(L)$ via formulas (\ref{fscales}).  Obviously, $\Upsilon(f_{q,a})(p,\gidoia)\leq h(p,\gidoia)$ whenever $p<0$ or $p\geq q$. If $0\leq p<q$, one has
 \[
  \Upsilon(f_{q,a})(p,\gidoia)=\left(\tbigvee^L_{r\geq p}\sigma_{q,a}(r) \right)^{\ast\ast}\leq h(q,\gidoia)^{\ast\ast}\leq h(p,\gidoia).
 \]
 Therefore $\Upsilon(f_{q,a})\in\mathcal{F}$. 
 
  Finally, by complete regularity and (r3), we have
 \[
 \begin{aligned}
 f_\vee(q,\gidoia)&\geq \tbigvee^L_{q^\prime>q}\tbigvee^L_{a\cbelow h(q^\prime,\gidoia)}\Upsilon(f_{q,a}(q,\gidoia))\\
 &\geq \tbigvee^L_{q^\prime>q}\tbigvee_{a\cbelow h(q^\prime,\gidoia)}a\\
 &=\tbigvee^L_{q^\prime>q}h(q^\prime,\gidoia)\\
 &=h(q,\gidoia).
 \end{aligned}
 \]
Hence, $f_\vee=h$.
 \end{proof}
We would like to emphasize the simplicity of this construction and the fact that it has no counterpart in the classical setting, as the Booleanization of a frame, even of a spatial one, is usually non-spatial.

  \section{An alternative representation: nearly finite Hausdorff continuous functions}
  
In order to present an alternative representation of the universal completion of $\mathrm{C}(L)$, we introduce now nearly finite Hausdoff continuous functions and show that they form a Riesz space in which the space of continuous real functions embeds naturally.
 
The set $\overline{\mathbb{IR}}$ of non-empty compact intervals $\boldsymbol{a}=[\underline{a},\overline{a}]$ of the extended real line ordered by reverse inclusion
 \[\boldsymbol{a}\sqsubseteq \boldsymbol{b}\quad\text{iff}\quad [\underline{a},\overline{a}]\supseteq [\underline{b},\overline{b}]\quad\text{iff}\quad\underline{a}\leq\underline{b}\leq \overline{b}\leq\overline{a}\]
 is a domain referred to, following Escard\'o \cite{Escardo97}, as the \emph{extended partial real line}.
This is  the extended version of the interval-domain that Dana Scott proposed as a domain-theoretic model for the real numbers \cite{S72}.  For general notions and results concerning domains we refer to \cite{compendium:03}, but we summarize here all what we need for this paper.
 The way-below relation in $\overline{\mathbb{IR}}$ is given by
 \[
 \boldsymbol{a}\ll\boldsymbol{b}\quad\text{iff}\quad \underline{a}<\underline{b}\leq\overline{b}<\overline{b}
 \]
 and we will denote $\{\boldsymbol{b}\mid \boldsymbol{a}\ll \boldsymbol{b}\}$ by $\twoheaduparrow \boldsymbol{a}$. The Scott topology $\mathcal{O}\overline{\mathbb{IR}}$ on $(\overline{\mathbb{IR}},\sqsubseteq)$ has a countable basis formed by the sets $\twoheaduparrow\boldsymbol{a}$ with $\underline{a},\overline{a}\in\Q$.  We will denote by $\pi_1,\pi_2\colon\overline{\mathbb{IR}}\to \overline{\mathbb{R}}$ the projections given by $\pi_1(\boldsymbol{a})=\underline{a}$ and $\pi_2(\boldsymbol{a})=\overline{a}$ for each $\boldsymbol{a}\in\overline{\mathbb{IR}}$. As $\pi_1$ is lower semicontinuous and $\pi_2$ is upper semicontinuous, we have that, for each $f\in\mathrm{C}(X,\overline{\mathbb{IR}})$, $\underline{f}=\pi_1\circ f$ is lower semicontinuous, $\overline{f}=\pi_2\circ f$ is upper semicontinuous and $\underline{f}\leq\overline{f}$. The extended real line $\overline{\R}$ embeds into $\overline{\mathbb{IR}}$ via the map $e\colon\overline{\R}\to\overline{\mathbb{IR}}$ given by $e(x)=[x,x]$ for each $x\in\overline{\R}$. Considering this we will identify $\overline{\R}$ with $e[\overline{\R}]\subseteq\overline{\mathbb{IR}}$ and a real-valued function $f\colon X\to\overline{\R}$ with $e\circ f\colon X\to\overline{\mathbb{IR}}$.
 
In general $\mathrm{C}(L)$ fails to be Dedekind complete due to axiom (r2). For that reason, the \emph{frame of partial reals} $\mathfrak{L}(\mathbb{IR})$ was introduced in \cite{MGP14} by removing this relation from the list (see also \cite{MThesis}). Specifically, $\mathfrak{L}(\mathbb{IR})$ is the frame generated by generators $(p,\gidoia)$ and $(\gidoia, q)$ for $p,q\in\Q$ subject to relations (r1), (r3)--(r6).  Its extended version, the \emph{frame of extended partial reals} $\mathfrak{L}(\overline{\mathbb{IR}})$ was introduced in \cite{MThesis} by additionally dropping (r5) and (r6). This is the pointfree version of the extended partial real line \cite{Escardo97}. It was shown in \cite{MThesis} that the spectrum $\Sigma\mathfrak{L}(\overline{\mathbb{IR}})$ is homeomorphic to $\overline{\mathbb{IR}}$ endowed with the Scott topology.

Following the arguments in \cite{MGP14} for the non-extended case of $\mathfrak{L}(\mathbb{IR})$, one can easily show that $\mathfrak{L}(\overline{\mathbb{IR}})$ is isomorphic to $\mathcal{O}\overline{\mathbb{IR}}$ and, therefore, it is spatial. Accordingly, frame homomorphisms $\mathfrak{L}(\overline{\mathbb{IR}})\to L$ are called \emph{continuous extended partial real functions}. The set $\overline{\mathrm{IC}}(L)$ of continuous extended partial real functions on $L$ is partially ordered by $f\leq g$ if
\[
f(p,\gidoia)\leq g(p,\gidoia)\quad\text{and}\quad f(\gidoia, q)\geq g(\gidoia, q)
\]
for all $p,q\in \Q$. Note that, in contrast with what happens in the case of continuous real functions, both conditions are necessary in this case. The assignments
\[
(p,\gidoia)\mapsto(p,\gidoia)\quad\text{and}\quad (\gidoia,q)\mapsto(\gidoia,q)
\]
determine the \emph{basic homomorphism} $\iota\colon\mathfrak{L}(\overline{\mathbb{IR}})\to\mathfrak{L}(\mathbb{R})$. The functions $f\in \overline{\mathrm{IC}}(L)$ that factor through the canonical homomorphism $\iota$ are just those that turn the defining relations (r2), (r5) and (r6) into identities in $L$. In view of this, we will  keep the notation $\mathrm{C}(L)$ to denote the class inside $\overline{\mathrm{IC}}(L)$ of functions that turn (r2), (r5) and (r6) into identities in $L$.

A function $f\in\overline{\mathrm{IC}}(L)$ is said to be \emph{Hausdorff continuous} if
\[
f(p,\gidoia)^\ast\leq f(\gidoia,q)\text{ and }f(\gidoia,q)^\ast\leq f(p,\gidoia)\text{ for all } p<q \text{ in }\Q. 
\]
We will denote by $\overline{\mathrm{H}}(L)$ the family of all Hausdorff continuous functions on $L$ \cite{M17} (see also \cite{MGP14}).

\begin{examples} \label{excharfunctions} For each $a,b\in L$ such that $a\wedge b=0$ let $\chi_{a,b}$ denote the continuous partial real function given by
\[
\chi_{a,b}(r,\gidoia)=
 \begin{cases}0&\text{if }r\ge 1,\\a&\text{if }0\le r< 1,\\ 1&\text{if }r< 0, \end{cases}
\quad\text{ and }\quad\chi_{a,b}(\gidoia,s)=
 \begin{cases}1&\text{if }s>1,\\ b&\text{if }0<s\le1,\\ 0&\text{if }s\le 0, \end{cases}
\]
for each $r,s\in{{\mathbb{Q}}}$. Clearly, $\chi_{a,b}\in\overline{\mathrm{H}}(L)$ if and only if $a=b^\ast$ and $b=a^\ast$ and $\chi_{a,b}\in\mathrm{C}(L)$ if and only if $a\vee b=1$, i.e., if and only if $a$ is complemented with complement $b$.
\end{examples}
  
  Hausdorff continuous functions can be alternatively described as the maximal elements of $\overline{\mathrm{IC}}(L)$ with respect to the following partial order on $\overline{\mathrm{IC}}(L)$:
  \[
  f\sqsubseteq g\quad\text{if}\quad f(p,\gidoia)\leq g(p,\gidoia)\text{ and }f(\gidoia, q)\leq g(\gidoia,q)\text{ for all }p,q\in\Q.
  \]
  
  The following two lemmas are the extended versions of \cite[Propositions 3.8 and 3.9]{MGP14} which refer to the lattice of continuous partial real functions $\mathrm{IC}(L)=\Frm(\mathfrak{L}(\mathbb{IR}),L)$. Since the defining relations (r5) and (r6) are not used at all in those propositions, the proofs are still valid in our case.

  \begin{lemma}
  The following are equivalent for any $f\in\overline{\mathrm{IC}}(L)$.
  \begin{enumerate}[\rm (1)]
  \item $f(p,\gidoia)^\ast\leq f(\gidoia,q)$ if $p<q$ in $\Q$.
  \item $f(\gidoia,q)=g(\gidoia,q)$ for all $q\in\Q$ and $g\in\overline{\mathrm{IC}}(L)$ such that $f\sqsubseteq g$.
  \end{enumerate}
  \end{lemma}
  
%
%
%
    \begin{lemma}
  The following are equivalent for any $f\in\overline{\mathrm{IC}}(L)$.
  \begin{enumerate}[\rm (1)]
  \item $f(\gidoia,q)^\ast\leq f(p,\gidoia)$ if $p<q$ in $\Q$.
  \item $f(p,\gidoia)=g(p,\gidoia)$ for all $p\in\Q$ and $g\in\overline{\mathrm{IC}}(L)$ such that $f\sqsubseteq g$.
  \end{enumerate}
  \end{lemma}
  
  \begin{corollary}
  A  function $f\in\overline{\mathrm{IC}}(L)$ is Hausdorff continuous if and only if it is maximal with respect to $\sqsubseteq$.
  \end{corollary}
  
  \begin{remark}\label{Hausdorffbikoitza}
  Let $f\in\overline{\mathrm{IC}}(L)$ be Hausdorff continuous and $r>p$ in $\Q$. Then
  \[
  f(r,\gidoia)^{\ast\ast}\leq f(\gidoia, r)^\ast\leq f(p,\gidoia).
  \]
  In consequence,
  \[
  f(p,\gidoia)=\tbigvee_{r>p}f(r,\gidoia)^{\ast\ast}
  \]
  for each $p\in\Q$,  like in the case of $\mathrm{C}(L)$ (Remark \ref{funtziobikoitza}).
  \end{remark}
  
 \begin{definition}
 We will say that $f\in\overline{\mathrm{H}}(L)$ is \emph{nearly finite} if $\tbigvee_{p\in\Q}f(p,\gidoia)$ and $\tbigvee_{q\in\Q}f(\gidoia,q)$ are dense elements in $L$, that is, if
 \[
\left(\tbigvee_{p\in\Q}f(p,\gidoia)\right)^{\ast\ast}=1=\left(\tbigvee_{q\in\Q}f(\gidoia,q)\right)^{\ast\ast}.
 \]
 We will denote by $\mathrm{H}_{nf}(L)$ the collection of all nearly finite functions on $L$.
 \end{definition}
 
 In \cite{M17} it was shown that $\overline{\mathrm{H}}(L)$ and $\overline{\mathrm{C}}(\mathfrak{B}(L))$ are order isomorphic with isomorphisms $\Gamma\colon\overline{\mathrm{H}}(L)\to\overline{\mathrm{C}}(\mathfrak{B}(L))$ where $\Gamma(f)=\beta_L\cdot f$ and $\Delta\colon\overline{\mathrm{C}}(\mathfrak{B}(L))\to \overline{\mathrm{H}}(L)$ where
\[
\Delta(g)(p,\gidoia)=\tbigvee^L_{r>p}g(r,\gidoia)\quad\text{and}\quad\Delta(g)(\gidoia,q)=\tbigvee^L_{s<q} g(\gidoia, s)
\]
for each $g\in\overline{\mathrm{C}}(\mathfrak{B}(L))$ and $p,q\in\Q$. As noted in \cite{M17}, these isomorphism  restrict to an isomorphism between $\mathrm{H}_{nf}(L)$ and $\mathrm{C}(\mathfrak{B}(L))$. Indeed, one has that $\tbigvee_{p\in\Q}f(p,\gidoia)$ and $\tbigvee_{q\in\Q}f(\gidoia,q)$ are dense elements in $L$ if and only if
  \[
\tbigvee^{\mathfrak{B}(L)}_{p\in\Q}f(p,\gidoia)=1=\tbigvee^{\mathfrak{B}(L)}_{q\in\Q}f(\gidoia,q).
 \]
 In consequence, $\mathrm{H}_{nf}(L)$ is a universally complete and an Archimedean Riesz space with respect to the operations inherited from $\mathrm{C}(\mathfrak{B}(L))$. Specifically, one has
 \[
 \lambda\odot f=\Delta(\lambda\cdot\Gamma(f))\quad\text{and}\quad f\oplus g=\Delta(\Gamma(f)+\Gamma(g)),\quad 
 \]
 for all $f,g\in\mathrm{H}_{nf}(L)$ and $\lambda\in\Q$. Naturally, the inverse of $f\in\mathrm{H}_{nf}(L)$, denoted by $-f$, is given by $\Delta(\Gamma(-f))$. It is straightforward to check the following.
 
 \begin{lemma}\label{operationsHnf}
  Let $0<\lambda\in\Q$ and $f,g\in\mathrm{H}_{nf}(L)$. For any $p,q\in\Q$, one has
  \begin{enumerate}[\rm(1)]
  \item $ (\lambda\odot f)(p,\gidoia)=f(\tfrac{p}{\lambda},\gidoia)$ and $(\lambda\odot f)(\gidoia, q)=f(\gidoia,\tfrac{q}{\lambda})$ for each $p,q\in\Q$.\smallskip
  \item  $(f\oplus g)(p,\gidoia)=\tbigvee_{r>p}\tbigvee_{t\in\Q}(f(t,\gidoia)\wedge g(r-t,\gidoia))^{\ast\ast}$ and\newline
$ (f\oplus g)(\gidoia, q)=\tbigvee_{s<q}\tbigvee_{t\in\Q}(f(\gidoia, t)\wedge g(\gidoia, s-t))^{\ast\ast}$.\smallskip
  \item $(-f)(p,\gidoia)=f(\gidoia,-p)$ and $(-f)(\gidoia, q)=f(-q,\gidoia)$.
  \end{enumerate}
 \end{lemma}
 
  \begin{corollary}\label{universalcompletion}
 For a completely regular frame $L$,  $\mathrm{H}_{nf}(L)$ is the universal completion of $\mathrm{C}(L)$.
 \end{corollary}

 \section{Universal completeness and isomorphic universal completions}
 
 For a complete regular frame $L$, 
 Banaschewski and Hong proved in \cite{BH03} (see also \cite{MGP14}) that $\mathrm{C}(L)$ is Dedekind complete iff $L$ is extremally disconnected and
Ball, Walters-Wayland and Zenk proved in \cite{BWZ11} that $\mathrm{C}(L)$ is $\sigma$-complete, i.e., every countable pairwise disjoint subset of $\mathrm{C}(L)^+$ has a supremum, iff $L$ is a $P$-frame. We provide here a characterization of completely regular frames $L$ for which $\mathrm{C}(L)$ is universally complete. 

Recall that an element $a\in L$ is said to be a \emph{cozero} if there exists an $f\in\mathrm{C}(L)$ such that $a=f(\gidoia,0)\vee f(0,\gidoia)$. In that case, there exists $g\in\mathrm{C}(L)$ such that $a=g(\gidoia,1)$ (simply take $g=((-f)\wedge f)+\boldsymbol{1}$). See \cite{BB97} for more details. A frame $L$ is said to be a \emph{$P$-frame} if each cozero element is complemented and said to be an \emph{almost $P$-frame} if $a=a^{\ast\ast}$ for each cozero element. Under extremal disconnectedness these two classes of frames coincide. Furthermore, $L$ is an almost $P$-frame if and only if $1$ is the only dense cozero \cite{Dube2009}.
 
 \begin{proposition}\label{isosuniversal}
 Let $L$ be a completely regular frame. TFAE:
 \begin{enumerate}[\rm (1)]
 \item $L$ is an extremally disconnected $P$-frame.
 \item $\mathrm{C}(L)=\mathrm{H}_{nf}(L)$.
 \item $\mathrm{C}(L)\simeq \mathrm{C}(\mathfrak{B}(L))$ as Riesz spaces.
 \item $\mathrm{C}(L)$ is universally complete. 
 \end{enumerate}
 \end{proposition}
 
 \begin{proof}
 \noindent (1)$\implies$(2): Let $f\in\mathrm{H}_{nf}(L)$ and $p<r<q$ in $\Q$. Then one has
 \[
 f(p,\gidoia)\vee f(\gidoia,q)\geq  f(r,\gidoia)^{\ast\ast}\vee f(r,\gidoia)^\ast=1
 \]
 if $L$ is extremally disconnected. Therefore $f\in\overline{\mathrm{C}}(L)$. Then $\tbigvee_{p\in\Q}f(p,\gidoia)$ and $\tbigvee_{q\in\Q}f(\gidoia,q)$ are cozeros which are by hypothesis dense. Therefore
 \[
 \tbigvee_{p\in\Q}f(p,\gidoia)=1=\tbigvee_{q\in\Q}f(\gidoia,q)
 \]
 if $L$ is a $P$-frame.
 
 \noindent (2)$\implies$(1):  For each $a\in L$, we have $\chi_{a^\ast,a^{\ast\ast}}\in\mathrm{H}_{nf}(L)$. Hence, if $\mathrm{C}(L)=\mathrm{H}_{nf}(L)$, 
 \[
 1=\chi_{a^\ast,a^{\ast\ast}}(0,\gidoia)\vee\chi_{a^\ast,a^{\ast\ast}}(\gidoia,1)=a^\ast\vee a^{\ast\ast}.
 \]
 Therefore $L$ is extremaly disconnected.
 
 Let $a$ be a dense cozero element of $L$. There exists an $f\in\mathrm{C}(L)$ such that $\boldsymbol{0}\leq f\leq\boldsymbol{1}$ and $a=f(\gidoia,1)$. Let $g\colon\mathfrak{L}(\overline{\mathbb{IR}})\to L$ be a frame homomorphism determined on generators by
 \[
  g(p,\gidoia)=f(\alpha(p),\gidoia)\quad\text{and}\quad g(\gidoia,q)=f(\gidoia, \alpha(q))
 \]
 for each $p,q\in\Q$ where $\alpha$ is an order isomorphism $\Q\to\Q\cap(-1,1)$. It is straightforward that such an assignment turns the defining relations (r1), (r3) and (r4) into identities in $L$ and that $g\in\overline{\mathrm{H}}(L)$. One has that
 \[
 \tbigvee_{p\in\Q}g(p,\gidoia)=\tbigvee_{p\in\Q}f(\alpha(p),\gidoia)\geq f\left(-\tfrac{1}{2},\gidoia\right)=1
 \]
 and
 \[
 \tbigvee_{q\in\Q}g(\gidoia, q)=\tbigvee_{q\in\Q}f(\gidoia, \alpha(q))=\tbigvee_{-1<q<1}f(\gidoia,q)=f(\gidoia,1)=a.
 \]
 Since $a$ is dense, $f\in\mathrm{H}_{nf}(L)$. If $\mathrm{C}(L)=\mathrm{H}_{nf}(L)$, we conclude that $a=1$ by $(r6)$. Therefore $L$ is a $P$-frame.
 
 \noindent (2)$\implies$(3): Obvious.
 
 \noindent (3)$\implies$(4): This follows from Corollary \ref{Booleanunivcomp}.
 
 \noindent (4)$\implies$(2): This follows directly from Theorem \ref{universalcompletion}.
 \end{proof}
 
 \begin{remarks}\label{remarkAlmostBoolean}
\noindent (1) With regard to the first condition in the proposition above, one may wonder whether all extremally disconnected $P$-frames are Boolean. This is not the case. In \cite{BB15} Banaschewski calls a frame $L$ \emph{almost Boolean} if it is regular, extremally disconnected and the Boolean algebra $B(L)$ of its completemented elements is closed under countable joins in $L$. He also shows that, for completely regular frames, almost Boolean frames are precisely the extremally disconnected $P$-frames. As the terminology suggests, the class of almost Boolean frames exceeds that of Boolean frames \cite{BB15}. Therefore, for  an almost Boolean frame $L$ which is not Boolean one has $\mathrm{C}(L)\simeq\mathrm{C}(\mathfrak{B}(L))$ while $L\not\simeq \mathfrak{B}(L)$.
In the spatial case, this issue reduces to a set theoretical question. See Section 6 below. 
\smallskip
 
\noindent (2) A ring $R$ is regular if for each $a\in R$ there exists a $b\in R$ such that $a=aba$. Banaschewski and Hong showed in \cite{BH03} that a frame $L$ is a $P$-frame iff $\mathrm{C}(L)$ is a regular ring. 
\end{remarks}

In summary, we obtain the following corollary.

\begin{corollary}
 Let $L$ be a completely regular frame. TFAE:
  \begin{enumerate}[\rm (1)]
 \item $L$ is almost Boolean.
 \item $\mathrm{C}(L)$ is Dedekind complete and regular as a ring.
 \item $\mathrm{C}(L)$ is universally complete. 
 \end{enumerate}
\end{corollary}
 
 We close this section with a direct consequence of the representations of the universal completion of $\mathrm{C}(L)$ presented in this paper.
 
 \begin{proposition}\label{twoframes}
 Let $L$ and $M$ be completely regular frames. TFAE:
 \begin{enumerate}[\rm (1)]
 \item $\mathfrak{B}(L)$ and $\mathfrak{B}(M)$ are isomorphic.
 \item $\mathrm{H}_{nf}(L)$ and $\mathrm{H}_{nf}(M)$ are Riesz  isomorphic.
 \item $\mathrm{C}(L)$ and $\mathrm{C}(M)$ have Riesz isomorphic universal completions.
 \end{enumerate}
 \end{proposition}
 
 \begin{proof}
 \noindent Since $\mathrm{C}(\mathfrak{B}(L))$ and $\mathrm{H}_{nf}(L)$ are isomorphic as Riesz spaces for all frames $L$, (1)$\implies$ (2) holds.  The reciprocal follows from Proposition \ref{isoBooleanpart}. (2)$\iff$(3) follows from Theorem \ref{universalcompletion}.
 \end{proof}

 \section{The spatial case}
Our pointfree approach sheds light also on the classical case of $\mathrm{C}(X)$ and provides a new construction of its universal completion. This construction is given in terms of interval-valued functions and arises from a direct lattice-theoretic approach.  In addition, we show how this new representation is equivalent to the one presented by van der Walt in \cite{W18}.
 
 By the (dual) adjunction between contravariant functors $\mathcal{O}\colon\Top\to\Frm$ and $\Sigma\colon\Frm\to\Top$, there is a natural isomorphism
 \[
 \Frm(L,\mathcal{O}X)\to\Top(X,\Sigma L)
 \]
 for any topological space $X$ and frame $L$. By the homeomorphism $\Sigma\mathfrak{L}(\overline{\mathbb{IR}})\simeq(\overline{\mathbb{IR}},\mathcal{O}\overline{\mathbb{IR}})$ (see \cite{MThesis}), we obtain that there is a natural isomorphism 
 \[
 \Psi\colon\overline{\mathrm{IC}}(\mathcal{O} X)\to \mathrm{C}(X,\overline{\mathbb{IR}})
 \]
given, for each $h\in\overline{\mathrm{IC}}(\mathcal{O}X)$, by
\[
\Psi(h)(x)=\left[
\tbigvee\{r\in\Q\mid x\in h(r,\gidoia)\},\tbigwedge\{r\in\Q\mid x\in h(\gidoia,r)\}\right]\quad\text{for all }x\in X.
\]

The following is a straightforward extension of \cite[Lemma 6.1]{MGP14}. Indeed, the proof in \cite{MGP14} for $\mathrm{IC}(\mathcal{O}X)$ is still valid for $\overline{\mathrm{IC}}(\mathcal{O}X)$ word by word, as defining relations (r5) and (r6) play no role there.

\begin{lemma} For any topological space $X$ and $f,g\in\overline{\mathrm{IC}}(\mathcal{O}X)$, one has:
\begin{enumerate}[\rm (1)]
\item $\pi_1\circ \Psi(f)\leq \pi_1\circ\Psi(g)$ if and only if $f(p,\gidoia)\leq g(p,\gidoia)$ for all $p\in\Q$.
\item $\pi_2\circ\Psi(f)\geq \pi_2\circ\Psi(g)$ if and only if $f(\gidoia, q)\leq g(\gidoia, q)$ for all $q\in\Q$.
\end{enumerate}
\end{lemma}

Consequently, as in the case of finite functions (see \cite{MGP14}), the map $\Psi$ is an order isomorphism for both $\leq$ and $\sqsubseteq$. Moreover, the restriction of $\Psi$ to $\mathrm{C}(\mathcal{O}X)\to\mathrm{C}(X)$ is also a Riesz space isomorphism. Just like in the pointfree case, we will say that a function $f\in\mathrm{C}(X,\overline{\mathbb{IR}})$ is \emph{Hausdorff continuous} if it is maximal with respect to $\sqsubseteq$. This definition is equivalent to the one given in \cite{RA04}. In addition, we will say that $f$ is \emph{nearly finite} if the set \[
f^{-1}\left(\twoheaduparrow[-\infty,+\infty]\right)=\{x\in X\mid \pi_1(f(x))\in\R\text{ and }\pi_2(f(x))\in\R\}
\]
is dense in $X$. Let us denote the set of all nearly finite Hausdorff continuous functions on $X$ by $\mathrm{H}_{nf}(X)$. The following fact follows immediately.

\begin{fact}
For any  $f\in\overline{\mathrm{IC}}(\mathcal{O}X)$ one has that
 \[f\in\mathrm{H}_{nf}(\mathcal{O}X)\quad\text{if and only if}\quad \Psi(f)\in\mathrm{H}_{nf}(X).
 \]
\end{fact}

 Therefore, $\mathrm{H}_{nf}(X)$ is a Riesz space with product with scalar $\odot$ and sum $\oplus$ operations induced by those from $\mathrm{H}_{nf}(\mathcal{O}X)$. Recall that a space $X$ is completely regular if and only if the frame $\mathcal{O}X$ is completely regular. The following follows immediately from Theorem \ref{universalcompletion}.

\begin{theorem}
For any completely regular space $X$, $\mathrm{H}_{nf}(X)$ is the universal completion of $\mathrm{C}(X)$.
\end{theorem}

In \cite{W18} van der Walt shows that,  for a completely regular space $X$, the Riesz space $\mathcal{NL}(X)$ of nearly finite normal lower semicontinuous functions is the universal completion of $\mathrm{C}(X)$. As the universal completion of a Riesz space is unique up to isormorphism, there has to be a Riesz isomorphism $\mathrm{H}_{nf}(X)\to\mathcal{NL}(X)$ that restricts to the identity homomorphism $\mathrm{C}(X)\to \mathrm{C}(X)$. In this respect, note that it follows directly from \cite[Theorem 1]{RA04}  that there is an order isomorphism between Hausdorff continuous functions $X\to\overline{\mathbb{IR}}$ and normal lower semicontinuous functions $X\to\overline{\mathbb{R}}$ given by $f\mapsto\underline{f}=\pi_1\circ f$. It is straightforward to check that this bijection restricts to a bijection between $\Pi\colon\mathrm{H}_{nf}(X)\to\mathcal{NL}(X)$. Unsurprisingly, this very bijection turns out to be the needed Riesz isomorphisms. It is immediate to see that $\Pi$ is an order isomorphism and that $\Pi(f)=f$ for every $f\in\mathrm{C}(X)$. Algebraic operations in $\mathcal{NL}(X)$ are given by 
\[
\lambda\odot f=\In\Su(\lambda\cdot f)\quad\text{and}\quad f\oplus g=\In\Su(f+g)
\]
for each $f,g\in\mathcal{NL}(X)$ and $\lambda\in\Q$, where $\cdot$ and $+$ are the usual pointwise operations on real functions and $\In$ resp.~$\Su$ denotes the lower resp.~upper Baire operator \cite{W18}.  The following technical lemma shows that $\Pi$ is indeed a Riesz space homomorphism.

 \begin{lemma}
  Let $0\leq \lambda\in\Q$ and $f,g\in\mathrm{H}_{nf}(X)$. Then
  \begin{enumerate}[\rm(1)]
  \item $\pi_1\circ(\lambda\odot f)=\In\Su(\lambda\cdot\underline{f}),$
  \item  $\pi_1\circ(f\oplus g)=\In\Su(\underline{f}+\underline{g}).$
  \end{enumerate}
  \end{lemma}
  
  \begin{proof}
 \noindent (1) Let $x\in X$. For any $\lambda>0$, one has  
 \[
 \begin{aligned}
 \pi_1((\lambda\odot f)(x))&=\pi_1(\Psi(\lambda\odot\Psi^{-1}(f))(x))\\
 &=\tbigvee\{r\in\Q\mid x\in(\lambda\odot f^{-1})(r,\gidoia)\}\\
 &=\tbigvee\{r\in\Q\mid x\in f^{-1}\left(\twoheaduparrow(\tfrac{r}{\lambda},+\infty]\right)\}\\
 &=\tbigvee\{r\in\Q\mid \underline{f}(x)> \tfrac{r}{\lambda}\}\\
 &=\lambda\underline{f}(x)=\In\Su(\lambda\cdot f)(x)
  \end{aligned}
 \]
 by Lemma \ref{operationsHnf}. For $\lambda=0$,  one has
 \[
 \pi_1((0\odot f)(x))=\pi_1(\Psi(0\cdot\Psi^{-1}(f))(x))=\pi_1(\Psi(\boldsymbol{0})(x))=0.
\]

%
%
 \noindent (2)  First note that
 \[
 \begin{aligned}
 \In\Su(\underline{f}+\underline{g})^{-1}(r,+\infty]&=\tbigcup_{p>r}\bigl((\underline{f}+\underline{g})^{-1}(p,+\infty]\bigr)^{-\circ}\\
 &=\tbigcup_{p>r}\left(\tbigcup_{t\in\Q}\underline{f}^{-1}(t,+\infty]\cap \underline{g}^{-1}(p-t,+\infty]\right)^{-\circ}
 \end{aligned}
 \]
 for any $r\in\Q$, where $A^-$ and $A^\circ$ denote the closure and interior of a subspace $A$ of $X$, respectively. For simplicity of notation, for each $r,t\in\Q$, let
 \[
a_{r,t}= \underline{f}^{-1}(t,+\infty]\cap\underline{g}^{-1}(p-t,+\infty].
\] 
 Then, for each $x\in X$,
 \[
 \begin{aligned}
 \In\Su(\underline{f}+\underline{g})(x)&=\tbigvee\left\{r\in\Q\mid x\in\tbigcup_{p>r}\left(\tbigcup_{t\in\Q}\underline{f}^{-1}(t,+\infty]\cap\underline{g}^{-1}(p-t,+\infty]\right)^{-\circ}\right\}\\
 &=\tbigvee\left\{r\in\Q\mid x\in\tbigcup_{p>r}\left(\tbigcup_{t\in\Q}a_{p,t}\right)^{-\circ}\right\}.
 \end{aligned}
 \]
 By  Lemma \ref{operationsHnf} and the fact that for an open subset $U$ of $X$, its double pseudocomplement $U^{\ast\ast}$ in $\mathcal{O}X$ is $U^{-\circ}$,
 \[
 \begin{aligned}
 \pi_1((f\oplus g)(x))&=\Psi(\Psi^{-1}(f)\oplus\Psi^{1}(g))(x)=\tbigvee\left\{r\in\Q\mid x\in(f^{-1}\oplus g^{-1})(r,\gidoia)\right\}\\
 &=\tbigvee\left\{r\in\Q\mid x\in\tbigcup_{p>r}\tbigcup_{t\in\Q}(f^{-1}(\twoheaduparrow(t,+\infty])\cap g^{-1}(\twoheaduparrow(p-t,+\infty]))^{-\circ}\right\}\\
 &=\tbigvee\left\{r\in\Q\mid x\in\tbigcup_{p>r}\tbigcup_{t\in\Q}(\underline{f}^{-1}(t,+\infty]\cap\underline{g}^{-1}(p-t,+\infty])^{-\circ}\right\}\\
 &=\tbigvee\left\{r\in\Q\mid x\in\tbigcup_{p>r}\tbigcup_{t\in\Q}a_{p,t}^{-\circ}\right\}
 \end{aligned}
 \]
 for each $x\in X$.
 
 Note that $ \pi_1((f\oplus g)(x))\leq \In\Su(\underline{f}+\underline{g})(x)$ follows from the fact that
\[
\left(\tbigcup_{t\in\Q}a_{p,t}\right)^{-\circ}\geq \tbigcup_{t\in\Q}a_{p,t}^{-\circ}.
\]
 For the inequality in the reverse direction, note that by Remark \ref{Hausdorffbikoitza}, for each $q>r$,
 \[
 (f^{-1}\oplus g^{-1})(q,\gidoia)^{\ast\ast}\leq  (f^{-1}\oplus g^{-1})(r,\gidoia),
 \]
that is,
\[
\left(\tbigcup_{p>q}\tbigcup_{t\in\Q}a_p,t^{-\circ}\right)^{-\circ}\subseteq\tbigcup_{p>r}\tbigcup_{r\in\Q}a_{p,t}^{-\circ}.
\] Therefore, if $r<q<s$, 
\[
\left(\tbigcup_{t\in\Q}a_{s,t}\right)^{-\circ}=\left(\tbigcup_{t\in\Q}a_{s,t}^{-\circ}\right)^{-\circ}
\subseteq\left(\tbigcup_{p>q}\tbigcup_{t\in\Q}a_p,t^{-\circ}\right)^{-\circ}\subseteq\tbigcup_{p>r}\tbigcup_{r\in\Q}a_{p,t}^{-\circ}.
\]
In consequence, $\In\Su(\underline{f}+\underline{g})(x)\leq \pi_1((f\oplus g)(x))$.
 \end{proof}
 
For the sake of completeness, we conclude this section with the application of Propositions \ref{isosuniversal} and \ref{twoframes} to the spatial case. For the latter, note that, in addition to slightly augmenting \cite[Corollary 39]{W18}, our proof avoids the use of the representation theory of Riesz spaces, and accordingly it is more direct and natural.
 
  \begin{corollary}
 Let $X$ and $Y$ be completely regular spaces. TFAE:
 \begin{enumerate}[\rm (1)]
 \item $\mathfrak{B}(\O X)$ and $\mathfrak{B}(\O Y)$ are isomorphic.
 \item $\mathrm{H}_{nf}(X)$ and $\mathrm{H}_{nf}(Y)$ are Riesz  isomorphic.
 \item $\mathcal{NL}(X)$ and $\mathcal{NL}(Y)$ are Riesz  isomorphic.
 \item $\mathrm{C}(X)$ and $\mathrm{C}(Y)$ have Riesz isomorphic universal completions.
 \end{enumerate}
 \end{corollary}

As each of the adjectives in the first condition of Proposition \ref{isosuniversal} is “conservative”, the result still holds if we replace frames by spaces. 
 \begin{corollary}
 Let $X$ be a completely regular space. TFAE:
 \begin{enumerate}[\rm (1)]
 \item $X$ is a extremally disconnected $P$-space.
 \item $\mathrm{C}(X)=\mathrm{H}_{nf}(X)$.
 \item $\mathrm{C}(X)=\mathcal{NL}(X)$.
 \item $\mathrm{C}(X)$ is universally complete. 
 \end{enumerate}
 \end{corollary}
 
The question whether extremally disconnected $P$-spaces are discrete was settled by Isbell in \cite{I55}: any extremally disconnected $P$-space of non-measurable power is discrete. However, the nonexistence of measurable cardinals is compatible with ZFC \cite{E96}.

\section{A pointfree approach to Maeda-Ogasawara-Vulikh representation theorem}

In this section we provide a localic version of the classical Maeda-Ogasawara-Vulikh representation theorem of Archimedean Riesz spaces with weak unit due to Maeda and Ogasawara \cite{MO42} and to Vulikh \cite{V47} independently. We recall here the main details of this classical result.  Let $R$ be an Archimedean Riesz space and let $\Omega$ denote the Stone space of its Boolean algebra of bands $\mathcal{B}(R)$. We call a continuous function $f\colon \Omega\to \overline{\R}$ nearly finite if it maps a dense subset of $\Omega$ into $\R$ and denote by $\mathrm{C}_{nf}(\Omega)$ the collection of all nearly finite extended continuous real functions on $\Omega$. Maeda-Ogasawara-Vulikh representation theorem states that $R$ embeds into $\mathrm{C}_{nf}(\Omega)$ and that this embedding constitutes its universal completion. A generalization to Archimedean lattice groups and Archimedean lattice rings was provided by Bernau in \cite{B65}. A localic version was given by Madden \cite{M90} (in fact, this is an application of a more general result for Archimedean $l$-groups with weak unit which can be found in \cite{M92, MV90}). This pointfree approach shows that Archimedean Riesz spaces with weak unit can be represented by continuous real functions on a frame. No extended valued functions are needed for the localic representation.  See \cite{F94} for a detailed survey of representation theorems in the theory of Archimedean Riesz spaces.

The aim of this section is to show that a  pointfree approach to Maeda-Ogasawara-Vulikh representation theorem yields a simpler and more natural proof than the classical one. For this purpose, we follow the proof as presented in \cite{LZ71} and take advantage of the following fact. $\mathcal{B}(R)$ is a frame, therefore a perfectly admissible ``space'' in the pointfree theory. Hence we do not need to construct its Stone space and this shortcut leads to clearer and more direct arguments. In contrast with \cite{M90}, where $R$ is represented by continuous real functions on the frame of relatively uniform closed ideals of $R$, we provide a representation by continuous real functions of its frame of bands.

For the convenience of the reader, we briefly summarize some definitions and facts concerning bands of Archimedean Riesz spaces. Let $R$ be an Archimedean Riesz space. An \emph{ideal} of $R$  is a Riesz subspace $A$ such that $|g|\leq |f|$ with $f\in A$ implies that $g\in A$. A \emph{band} of $R$ is an ideal $A$ such that whenever a subset of $A$ has a supremum in $R$, then that supremum is in $A$. We will denote by $\mathcal{B}(R)$ the set of all bands of $R$. Ordered by inclusion, $\mathcal{B}(R)$ forms a complete Boolean algebra with infima given by  intersection. For $f\in R$ (resp. $A\subseteq R$), we will denote by $[f]$ (resp. $[A]$) the band generated by $f$ (resp. $A$), that is, the smallest band containing $f$ (resp. $A$). In particular one has
that $g\in[f]$ if and only there exists $\{f_n\}_{n\in\N}$ such that $f_n\leq n|f|$ for each $n\in\N$ and 
\[
|g|\leq \tbigvee_{n\in\N}f_n.
\]
 Further, for each $f,g\in R^+$, one has:
\begin{enumerate}[(1)]
\item $[f]\vee[g]=[f+g]=[f\vee g]$
\item $[f]\wedge[g]=[f\wedge g]$.
\end{enumerate}
For $0<f,g\in R$, one has that $g$ is in $[f]^\ast$ (the complement of $[f]$ in $\mathcal{B}(R)$) if $f\wedge g=0$. A \emph{weak unit} of $R$ is an element $0< e\in R$ such that the band generated by $e$ is $L$ itself. Equivalently, $e$ is a weak unit if $e\wedge f$ implies  $f=0$. For general results concerning Archimedean Riesz spaces we refer to Luxenburg-Zaneen \cite{LZ71}.

\begin{proposition}
Let $R$ be an Archimedean Riesz space with a weak unit $e$ and $f\in R$. The map $\sigma_f\colon\Q\to\mathcal{B}(L)$ given by
\[
p\mapsto [(f-pe)^+]
\]
is a scale.
\end{proposition}

\begin{proof}
As $\sigma_f$ is obviously antitone and every image is complemented, since $\mathcal{B}(L)$ is a Boolean algebra, $\sigma_f$ is an extended scale. We shall check that $\tbigvee_{p\in\Q}[(f-pe)^+]=1$. As $e$ is a weak unit, there exists $\{f_n\}_{n\in\Q}\subseteq R$ such that $0\leq f_n\leq ne$ and $|e-f|\leq\tbigvee_{n\in \N}f_n$. Consequently $ e-f\leq \tbigvee_{n\in\N}f_n$. Thus 
\[
e\leq f+\tbigvee_{n\in\N}f_n\leq\tbigvee_{n\in\N}(f+f_n).
\]
Accordingly, one has
\[
1=\tbigvee_{n\in\N}[(f+f_n)^+]\leq\tbigvee_{n\in\N} [(f+ne)^+]\leq\tbigvee_{p\in\Q}[(f-pe)^+].
\]
In order to check that $\tbigvee_{p\in\Q}[(f-pe)^+]^\ast=1$, simply note that, as  $(f-pe)^-\wedge(f-pe)^+=0$, one has that $(f-pe)^-\in[(f-pe)^+]^\ast$. Therefore
\[
\tbigvee_{p\in\Q}[(f-pe)^-]\leq\tbigvee_{p\in\Q}[(f-pe)^+]^\ast
\]
Since $(f-pe)^-=(-f+pe)^+$, we can use the same argument as before applied to $-f$ and show that $\tbigvee_{p\in\Q}[(f-pe)^-]=1$.
\end{proof}

Therefore, for each element $f\in R$, the scale $\sigma_f$ determines a continuous real function on the Boolean algebra of bands of $R$. We will denote this function by $m(f)$.
\begin{proposition} Let $R$ an Archimedean Riesz space with weak unit $e$.
The map $m\colon R\to \mathrm{C}(\mathcal{B}(R))$ given by $f\mapsto m(f)$ is a Riesz space embedding.
\end{proposition}

\begin{proof}
Let $\lambda,p\in\Q$ and $f,g\in R$. One has
\[
[(\lambda f-pe)^+]=\left[\left(f-\tfrac{p}{\lambda} e\right)^+\right].
\]
Therefore $\lambda m(f)=m(\lambda f)$.

For any $t\in\Q$,
\[
\begin{aligned}
[(f+g-pe)^+]&=[(f+g-pe+te-te)^+]\\
&=[((f-te)+(g-(p-t)e))^+]\\
&\geq [((f-te)\wedge (g-(p-t)e))^+]\\
&=[(f-te)^+]\wedge[(g-(p-t)e)^+].
\end{aligned}
\]
Thus $m(f)+m(g)\leq m(f+g)$. Applying this result to $-f$ and $-g$ we obtain that the opposite inequality holds. Therefore, $m(f)+m(g)=m(f+g).$

In order to check that $m$ is an embedding, let $f\in\R$ be such that $m(f)=0$. Then, for any $p>0$, one has $[(f-pe)^+]=0$. Thus $f\leq pe$ for any $p>0$. As $R$ is Archimedean, we conclude that $f\leq 0$. In addition, as obviously $(f-pe)^-\wedge(f-pe)^+=0$, one has that $[(f-pe)^-]=0$ for any $p<0$. Then, $-f\leq -pe$ for any $p<0$. Again, as $R$ is Archimedean, $-f\leq 0$. In conclusion $f=0$.
\end{proof}

\begin{proposition}
Let $R$ be an Archimedean Riesz space with weak unit $e$. Then the ideal $D$ generated in $\mathrm{C}(\mathcal{B}(R))$ by $m(R)$ is the Dedekind completion of $R$ and $[D]=\mathrm{C}(\mathcal{B}(R))$.
\end{proposition}

\begin{proof}
First note that 
\[
D=\{g\in\mathrm{C}(\mathcal{B}(R))\mid |g|\leq m(f)\text{ for some } f\in R\}.
\]
As $\mathrm{C}(\mathcal{B}(R))$ is Dedekind complete, $D$, as a Riesz space, is also Dedekind complete. In order to check that $D$ is the Dedekind completion of $R$ it is sufficient to check that for each $0<g\in D$ there exist $0\neq f,h\in L$ such that
\[
0\leq m(f) \leq g \leq m(h),
\]
by \cite[Theorem 32.6]{LZ71}.
The existence of $h$ follows from the definition of $D$. As $g\neq 0$, there exists $p\in\Q$ such that $g(p,\gidoia)\neq \{\boldsymbol{0}\}$. Thus there exists $0< f\in g(p,\gidoia)$. Without loss of generality, we can assume that $f\leq p e$. Then, for any $0\leq r<p$ in $\Q$,
\[
m(f)(r,\gidoia)\leq m(f)(0,\gidoia)\leq [f]\leq g(p,\gidoia)\leq g(r,\gidoia).
\]
As $m(f)\leq\boldsymbol{p}$, then $m(f)(r,\gidoia)=0$ for all $r\geq p$. Consequently, $0<m(f)\leq g$.

Finally $[D]=\mathrm{C}(\mathcal{B}(R))$ as the band generated by $m(e)=\mathbf{1}$ already generates $\mathrm{C}(\mathcal{B}(R))$.
\end{proof}

As a corollary we obtain the following. 

\begin{theorem}
Let $R$ be an Archimedean Riesz space with weak unit $e$. The embedding $m\colon R\to \mathrm{C}(\mathcal{B}(R))$ constitutes the universal completion of $R$.
\end{theorem}

Finally, by Proposition \ref{isoBooleanpart} and the uniqueness of the universal completion, we obtain the following result.

\begin{corollary}
For any completely regular frame $L$, the Booleanization $\mathfrak{B}(L)$ of $L$ is isomorphic to the frame of bands $\mathcal{B}(\mathrm{C}(L))$ of $\mathrm{C}(L)$.
\end{corollary}

\medskip
\noindent{\bf Acknowledgements.}
 The author wishes to thank Frederick Dashiell for his many helpful suggestions during the preparation of the paper. The author gratefully acknowledges support from the Spanish Ministry of Science, Innovation and Universities, reference code PID2019-103838GB-100 (MCIU/AEI/FEDER, UE) and from the Basque Government (Postdoctoral Fellowship POS\_2017\_2\_0042 and grant IT974-16).

\end{document}